\documentclass[12pt]{article}
\usepackage{amsmath}
\usepackage{amssymb}
\usepackage{verbatim}
\usepackage{graphicx}

\title{A TEST AGAINST TREND IN RANDOM SEQUENCES}
\author{Peter Lindqvist}
\date{\footnotesize{Department of Mathematical Sciences\\Norwegian University of Science and Technology\\NO--7491 Trondheim, Norway}}
\begin{document}
\maketitle

%\begin{center}
 % \texttt{Dedicated to Olli Martio on his seventy-fifth birthday}
%\end{center}
 \bigskip
{\small \textsc{Abstract:}}\footnote{AMS classification 62G60} \textsf{We study a modification of Kendall's $\tau$, replacing his permutations of $n$ different numbers by sequences of length $n$. Thus repetition is allowed. In particular, binary sequences are studied.}

\section{Introduction}\label{introduction}
The basic tool in Kendall's $\tau$-test is the ``score'' $S$. Suppose that $\ell$ digits $0,1,...,\ell -1$ satisfying the transitive relations $0 < 1 <...<\ell-1$ are given and consider all the $\ell^n$ possible sequences
\begin{equation}\label{x}
  x_1,x_2,...,x_n
  \end{equation}
of length $n$ that can be formed by the aid of these digits. Let $S^+(x_1,...,x_n)$ denote the number of true inequalities
$$x_{i} > x_{j}\quad\text{with}\quad i>j$$
in the sequence (\ref{x}). Analogously, $S^-(x_1,...,x_n)$ counts the number of all valid inequalities $x_i <x_j$ with $i>j$. Define
\begin{equation}\label{S}
  S(x_1,...,x_n) = S^+(x_1,...,x_n) - S^-(x_1,...,x_n).
\end{equation}
For example, $S = 5-11 = -6$ for the sequence $0\,1\,1\,2\,0\,2\,1.$
Originally\footnote{The coefficient $\tau$ was considered by Greiner (1909) and Esscher (1924). The coefficient was rediscovered by Kendall (1938).}, Kendall considered the distribution of $S$ among the $\ell !$ permutations of the digits $\,0,1,...,\ell-1\,$ (then $n = \ell$), and so far as I know\footnote{However, see \cite{L1}.} the generalizations of Kendall's $\tau$-test rely upon the distribution of $S$ among all
$$\dfrac{n!}{(2!)^{p_{2}}(3!)^{p_{3}}\cdots (r!)^{p_r}}\qquad (p_1+2p_2+...+rp_r = n)$$ possible sequences (\ref{x}) consisting of $p_1$ digits occurring only once, $p_2$ pairs, $p_3$ triplets, and so on. Here the numbers of ties, i.e., $p_1,p_2,...,p_r$, are regarded as fixed. See \cite{S}. For permutations $S^+$ has been thoroughly investigated in \cite{M}.

We shall study a different situation, arising for example in connexion with the testing of \emph{sequences of random digits}. In this setting the number of ties cannot be regarded as fixed  \emph{a priori}. Thus we are led to study \emph{the distribution of $S$ among all $\ell^n$ sequences (\ref{x})}. As we shall learn,\emph{ this distribution approaches normality, as $n \to \infty.$}

Having in mind applications for a certain kind of sampling, we have considered the binary case $\ell = 2$ also when the probability for a $0$ is $p$ and the probability for a $1$ is $q$, where $p+q=1.$ See Section $\ref{Binary}$.

Finally, we mention that the distribution for $S$ in our setting of the problem is, in certain respects even simpler than the version considered by Kendall \cite{K1}, Sillitto \cite{S}, and Silverstone \cite{S2}.

\section{The Basic Results}\label{Basic R}

The mean value for $S$ taken over all sequences (\ref{x}) is zero by symmetry:
\begin{equation}
  \mu(n)\,=\,0.
\end{equation}
When all the sequences are equiprobable, the variance is
\begin{equation}\label{var}
  \sigma^2(n)\,=\, \frac{\ell-1}{\ell}\,\frac{n(n-1)}{2}\,+\,\frac{\ell^2-1}{\ell^2}\,\frac{n(n-1)(n-2)}{9}
\end{equation}
and the fourth central moment is
\begin{align}\label{four}  \mu_4(n)\,&=\,\Bigl(\frac{\ell^2-1}{\ell^2}\Bigr)^2\frac{100n^4+328n^3-127n^2-997n-372}{2700}\,n(n-1)\nonumber\\  &\quad+\frac{\ell^2-1}{\ell^4}\,\frac{252n^3+507n^2-3623n+3652}{900}\,n(n-1)\nonumber\\
  &\quad-\frac{\ell^2-1}{\ell^3}\,\frac{2n^3+3n^2-5n-15}{6}\,n(n-1)\\
  &\quad+\frac{\ell-1}{\ell^3}\,\frac{n^2+11n-25}{2}\,n(n-1).\nonumber
\end{align}
By symmetry
$$\mu_3(n)\,=\,0,\, \mu_5(n)\,=\,0,\, \mu_7(n)\,=\,0,....$$
It is interesting to observe that for $n$ fixed the moments approach those given by Kendall in \cite{K1}, as $\ell \to \infty.$ Formula (\ref{var}) is derived in Section \ref{The Vaiance}, but the corresponding calculations for (\ref{four}) are, to say the least, a laborious task and so $\mu_4(n)$ is given without proof, when $\ell \geq 3.$

As $n$ grows, the distribution for $S$ tends towards normality in the sense that the frequency between the values $S_1$ and $S_2$ tends to
$$\frac{1}{\sigma(n)\sqrt{2\pi}}\int_{S_1}^{S_2}\!e^{-x^2/2\sigma^2(n)}\,dx,$$
  where the standard deviation is $\sigma(n)=\sqrt{\mu_2(n)}.$ This follows from the Second Limit Theorem, since
  \begin{equation}\label{normal}
    \lim_{n\to\infty}\frac{\mu_{2k}(n)}{(\sigma(n))^{2k}}\,=\, 1\cdot2\cdot5\cdots (2k-1)
  \end{equation}
  and the odd moments are zero. It is easy to prove (\ref{normal}) for small values of $\ell$, but the probability function for $S$ becomes soon too complicated, as $\ell$ grows. Therefore we shall prove (\ref{normal}) only for the binary case $\ell = 2$, see Section \ref{approach}.

  In the binary case the probability generating function for $S$ is
  \begin{equation}\label{f}   f(x)\,=\,
    \begin{cases}
      \frac{1}{2^{2\nu}}\,\underset {k=1}{\overset{k=\nu}{\prod}}\bigl(x^{2k-1}+2+x^{1-2k}\bigr),\qquad n = 2\nu\\
       \frac{2}{2^{2\nu+1}}\,\underset{k=1}{\overset{k=\nu}{\prod}}\bigl(x^k+x^{-k}\bigr)^2,\qquad\qquad \,\,\,  n = 2\nu+1
  \end{cases}\end{equation}
  and so the characteristic function $\phi(\theta)\,=\,f(e^{i\theta})$ reduces to the simple expression (\ref{char}). Our proof for (\ref{normal}), when $\ell=2$ is based on $\phi(\theta)$. Furthermore, the distribution for $S$ can be rapidly calculated via a suitable interpretation of (\ref{f}).

  In the general binary case, when the probability that $x_i=0$ is $p$ and that $x_i =1$ is $q$ in (\ref{x}), $i = 1,2,...,n,\,\,p+q =1,$ the corresponding probability function is given by (\ref{pqgen}) and the characteristic function by (\ref{pqchar}). Now again $\mu(n) =0$, and
  \begin{equation}\label{varpq}
    \sigma^2(n)\,=\,\frac{n(n^2-1)}{3}\,pq
  \end{equation}
  and the fourth moment is
  \begin{align}\label{fourpq}
    \mu_4(n)\,&=\, \frac{n(n^2-1)(5n^3-6n^2-5n+14)}{15}\,p^2q^2\nonumber \\
    &\quad+\,\frac{n(n^2-1)(3n^2-7)}{15}\,pq(p-q)^2.
  \end{align}
  All odd moments are zero and the asymptotic normality  (\ref{normal}) holds even for $p\neq q$.

  \section{The Probability Function}\label{probability}

  Consider the binary case $\ell = 2$ with equiprobable sequences (\ref{x}). Direct calculation of $S$ for $n=3$ yields
  \begin{eqnarray*}
    0\,0\,0\,\,\mathbf{\phantom{-}0}\quad 0\,1\,0\,\,\mathbf{\phantom{-}0}\quad1\,0\,0\,\,\mathbf{2}\quad1\,1\,0\,\,\mathbf{2}\\ 0\,0\,1\,\,\mathbf{-2}\quad 0\,1\,1\,\,\mathbf{-2}\quad1\,0\,1\,\,\mathbf{0}\quad 1\,1\,1\,\,\mathbf{0}
  \end{eqnarray*}
  and so we can write
  \begin{eqnarray*}
  \underline{\mathbf{-2-\!1\,\,\,0\,+\!1\,+\!2}}\\
  \phantom{-}2\,\,\phantom{-}0\,\,\,4\quad\,\,0\quad\,\, 2
  \end{eqnarray*}
  and in this manner the distribution for $S$ can be tabulated. The results are displayed below:
  \begin{gather*}
    2\\1\,2\,1\\2\,0\,4\,0\,2\\1\,2\,1\,2\,4\,\,2\,1\,2\,1\\
    2\,0\,4\,0\,6\,0\,8\,0\,6\,0\,4\,0\,2\\
    1\,2\,1\,2\,4\,4\,5\,4\,5\,8\,5\,4\,5\,4\,4\,2\,1\,2\,1
  \end{gather*}
  The fundamental observation is that the table can be constructed via the following kind of figurates. For even $n$ we have
  \begin{gather*}
    \overline{1\,2\,1} \\
    1\,2\,1\\
   1\,2\,1\,\phantom{1\,2\,1}\,1\,2\,1\\
   \overline{ 1\,2\,1\,2\,4\,2\,1\,2\,1}
    \\
    1\,2\,1\,2\,4\,2\,1\,2\,1\\
   1\,2\,1\,2\,4\,2\,1\,2\,1\,\phantom{4}\, 1\,2\,1\,2\,4\,2\,1\,2\,1\\
   \overline{ 1\,2\,1\,2\,4\,4\,5\,4\,5\,8\, 5\,4\,5\,4\,4\,2\,1\,2\,1}
  \end{gather*}
  and so on. (The overlined sequences display $n=2,\,4$ and $6$.) For odd $n$ the table looks like
  
     \begin{gather*}
    \overline{0\,2\,0} \\
    0\,2\,0\\
   0\,2\,0\,\phantom{1}\,0\,2\,0\\
   \overline{ 0\,2\,0\,4\,0\,2\,0}
    \\
    0\,2\,0\,\,4\,0\,2\,0\\
   0\,2\,0\,4\,0\,2\,0\,\phantom{4}\, 0\,2\,0\,4\,0\,2\,0\\
   \overline{ 0\,2\,0\,4\,0\,6\,0\,8\,0\, 6\,0\,4\,0\,2\,0}\\
   0\,2\,0\,4\,0\,6\,0\,8\,0\, 6\,0\,4\,0\,2\,0\\
   0\,2\,0\,4\,0\,6\,0\,8\,0\, 6\,0\,4\,0\,\mathbf{4}\,0\,4\,0\,6\,0\,8\,0\, 6\,0\,4\,0\,2\,0\\
   \overline{\phantom{2\,0\,6\,0\,12\,0\,14\,0\,16\,0\,20\,0\,16\,14\,12\,0\,6\,0\,4\,0\,2}}
     \end{gather*}
     and the next row, corresponding to $n=7$, becomes
     $$
     0\, 2\,0\,4\,0\,6\,0\,12\,0\,14\,0\,16\,0\,20\,0\,16\,0\,14\,0\,12\,0\,6\,0\,4\,0\,2\,0.$$
    (The first and last zeros in a row are void.) In an obvious interpretation the above process reads
     \begin{align*}
       &x^{-1}+2+x,\\
       &x^{-4}+2x^{-3}+x^{-2}+2x^{-1}+4+2x+x^2+2x^3+x^4\\
      &\quad =\, \bigl(x^{-1}+2+x\bigr)\bigl(x^{-3}+2+x^3\bigr), .....
     \end{align*}
     for odd $n$ and
     \begin{align*}
      &2,\,\,2x^{-2}+4+2x^2,\\
       &2x^{-6}+4x^{-4}+6x^{-2}+8+6x^2+4x^4+2x^6\\
       &\quad =2\bigl(x^{-2}+2+x^2\bigr)\bigl(x^{-4}+2+x^4\bigr), .....
     \end{align*}
     for even $n$. This leads to the probability generating function (\ref{f}).

     A simple proof for the probability generating function $f(x)$ comes from considering the binary sequence
     $$j_1,j_2,...,j_n\qquad \text{where} \quad  j_k = 0\quad \text{or} \quad 1.$$
     Then we have
     $$S\,=\, \sum_{k=2}^n(j_1-j_k)\,+\sum_{k=3}^n(j_2-j_k)\,+\,...\,+\sum_{k=n}^n(j_{n-1}-j_k)$$
     and the index $j_k$ appears exactly $(n-k)-(k-1)$ times and so its contribution to the score $S$ is
     $$[(n-k)-(k-1)]j_k.$$
     Therefore
     $$S\,=\,(n-1)j_1+(n-3)j_2+\cdots+(n-2k+1)j_k+\cdots+(3-n)j_{n-1}+(1-n)j_n$$
     for this sequence. Now $j_k$ is $0$ or $1$ so that the generating function becomes     $$\bigl(1+x^{n-1}\bigr)\bigl(1+x^{n-3}\bigr)\cdots\bigl(1+x^{-(n-3)}\bigr)\bigl(1+x^{-(n-1)}\bigr).$$
     Upon multiplication, the coefficient of $x^t$ indicates how many times $S=t$ among all possible sequences $j_1,j_2,...,j_n$. Dividing by the total number of sequences we arrive at the probability generating function $(\ref{f})$.

     The characteristic function for $S$ is $\phi(\theta)\,=\,f(e^{i\theta}),\,\,i^2 = -1.$ Euler's formula yields
     \begin{equation}\label{char}
       \phi(\theta)\,=\begin{cases}
       \underset{k=1}{\overset{\nu}{\prod}}\cos^2(k\theta),\qquad n=2\nu+1\\
       
       \underset{k=1}{\overset{\nu}{\prod}}\cos^2\bigl(\frac{2k-1}{2}\,\theta\bigr),\quad
       n=2\nu.
       \end{cases}
     \end{equation}
     By definition $\phi(0) =1$ and
     $$\phi^{(k)}(0)\,=\,i^{-k}\mu_k(n),\qquad k=1,2,...,n.$$
     By symmetry $\phi'(0)=0,\,\phi^{3}(0)=0$,..., so that all odd moments are zero and
     $$\phi^{(2k)}(0)\,=\,(-1)^k\mu_{2k}(n).$$
     Direct calculations yield
     \begin{align*}
       &\mu_1(n)\,=\,0,\\
       &\mu_2(n)\,=\,\frac{n(n^2\!-\!1)}{12},\\
       &\mu_3(n)\,=\,0,\\
       &\mu_4(n)\,=\,\frac{n(n^2\!-\!1)(5n^3\!-\!6n^2\!-\!5n\!+\!14)}{240}\\
       &\mu_5(n)\,=\,0,\\
       &\mu_6(n)\,=\frac{n(n^2\!-\!1)(35n^6\!-\!126n^5\!+\!74n^4\!+\!420n^3\!-\!829n^2\!-\!294n\!+\!1488)}{4032},\\
       &\mu_7(n)\,=\,0,\\
       &.......
     \end{align*}
 The arrangements in Section \ref{approach}  will shorten such calculations.
     
     \section{Approach to Normality}\label{approach}

     In order to show that the distribution for $S$ approaches normality in the binary case, we shall prove (\ref{normal}). The dichotomy in formulae (\ref{char}) forces us to separate the cases
     $$\lim_{\nu \to \infty}\frac{\mu_{2k}(2\nu)}{\bigl(\sigma(2\nu)\bigr)^{2k}}\,=\, \frac{(2k)!}{2^k\cdot k!}, \qquad
     \lim_{\nu \to \infty}\frac{\mu_{2k}(2\nu+1)}{\bigl(\sigma(2\nu+1)\bigr)^{2k}}\,=\, 
       \frac{(2k)!}{2^k\cdot k!}.$$
       However, both cases are so similar that we shall write down only the odd case $n= 2\nu+1.$ Then the characteristic function is
       $$\phi(\theta)\,=\,cos^2(\theta)\, cos^2(2\theta) \cdots  cos^2(\nu \theta) $$
       and by logarithmic differentiation
       \begin{equation}\label{tang}
         \phi'(\theta)\,=\,-2\,\phi(\theta)\sum_{j=1}^{\nu}j\,\tan(j\theta).
       \end{equation}

       Denoting
       \begin{equation}\label{A}
         A(\theta)\,=\,-2\sum_{j=1}^{\nu}j\,\tan(j\theta),
       \end{equation}
       we obviously have
       $$A(0)\,=\,0,\, A''(0)\,=\,0,\, A^{(4)}(0)\,=\,0,....$$
       In order to calculate the odd derivatives $A'(0),\, A'''(0),...$ we use the expansion
       $$ \tan(z)\,=\,\sum_{k=1}^{\infty}\frac{2^{2k}\bigl(2^{2k}-1\bigr)}{(2k)!}\,(-1)^{k+1}\mathrm{B}_{2k}z^{2k-1}\qquad \Bigl(|z|^2 < \frac{\pi^2}{4}\Bigr)$$
       where $\mathrm{B}_0 =1, \mathrm{B}_2 =1/6,   \mathrm{B}_4 =-1/30,... $
       are the Bernoulli numbers. Thus
       \begin{equation}\label{Aexp}
         A(\theta)\,=\,\sum_{k=1}^{\infty}\frac{2^{2k}\bigl(2^{2k}-1\bigr)}{(2k)!}\,(-1)^{k}\mathrm{B}_{2k}\bigl(1^{2k}+2^{2k}+\dots+\nu^{2k}\bigr)\theta^{2k-1}
       \end{equation}
       for $|\theta| < \pi/2\nu.$ We deduce that
       \begin{equation}\label{Azero}
         A^{(2k-1)}(0)\,=\,\frac{(-1)^k}{k}\,2^{2k}\bigl(2^{2k}-1\bigr)\mathrm{B}_{2k}\,s_{\nu}(2k)
       \end{equation}
       where
       \begin{align*}
         s_{\nu}(2k)\,& =\, 1^{2k}+2^{2k}+\dots+\nu^{2k}\\
         &=\,\frac{\nu^{2k+1}}{2k+1}+ \frac{\nu^{2k}}{2} +\frac{k\nu^{2k-1}}{6} + \langle\mathrm{lower\, terms}\rangle
       \end{align*}
       are well-known polynomials of degree $2k+1$.

       According to (\ref{tang}) and (\ref{A}) we obtain\footnote{The connexion with cumulants is obvious, since $\frac{d}{d\theta}\log|\phi(\theta)|= A(\theta).$} by Leibniz rule
       \begin{align*}
         \phi'(\theta)\,&=\,\phi(\theta)A(\theta)\\
          \phi''(\theta)\,&=\,\phi'(\theta)A(\theta)+\phi(\theta)A'(\theta)\\
         \phi'''(\theta)\,&=\,\phi''(\theta)A(\theta) +2\phi'(\theta)A'(\theta)
         +   \phi(\theta)A''(\theta)\\
         &\quad\vdots\\
         \phi^{k+1}(\theta)\,&=\,\sum_{j=0}^k\binom{k}{j}\phi^{(j)}(\theta)A^{(k-j)}(\theta)\\
         &\quad\vdots
       \end{align*}
       In passing, we calculate
       \begin{gather*}
         \phi''(0)\,=\,A'(0)\,=\,-\frac{n\bigl(n^2-1\bigr)}{12}\\
         \phi^{(4)}\,=\,3\phi''(0)^2+A'''(0)\,=\,\frac{n\bigl(n^2-1\bigr)\bigl(5n^3-6n^2-5n+14\bigr)}{240}
       \end{gather*} for odd $n$. For even $n$ we shall arrive at the same formulae. We have obtained that
       \begin{align}\label{sigma}
         \sigma^2(n)\,&=\, \frac{n\bigl(n^2-1\bigr)}{12}\\
         \mu_4(n)\,&=\,\frac{n\bigl(n^2-1\bigr)\bigl(5n^3-6n^2-5n+14\bigr)}{240}.
       \end{align}
       This shows that (\ref{normal}) holds at least for $k=1$ and $k=2$. For general $k$ we use induction.

       To this end, notice that at the point $\theta = 0$ we have
       $$\phi^{(2k)}\,=\,(2k-1)\phi''\phi^{(2k-2)}    +\sum_{j=2}^k\binom{2k-1}{2j-1}A^{(2j-1)}\phi^{(2k-2j)}$$
         or, more conveniently,
       \begin{equation}\label{induction}
       \frac{\phi^{(2k)}}{(\phi'')^k}\,=\,(2k-1)  \frac{\phi^{(2k-2)}}{(\phi'')^{k-1}}
           +\sum_{j=2}^k\binom{2k-1}{2j-1} \frac{\phi^{(2k-2j)}}{(\phi'')^{k-j}}\frac{ A^{(2j-1)}}{(\phi'')^j}.
       \end{equation}
             According to (\ref{Azero}) and (\ref{sigma}), where $n = 2\nu+1$ is odd, we have
             $$\lim_{\nu \to \infty}\frac{ A^{(2j-1)}}{(\phi'')^j}\,=\,0,\qquad j=2,3,...,k,$$
             since $A^{(2j-1)} \approx \nu^{2j+1}$ and $(\phi'')^j \approx \nu^{3j}.$
             But now (\ref{induction}) shows that, if
             $$\lim_{\nu\to\infty}\frac{\phi^{(2m)}}{(\phi'')^m}\,=\,1\cdot2\cdot\cdots (2m-1)$$
             holds for $m=1,2,...,k-1,$ then also
             $$\lim_{\nu\to\infty}\frac{\phi^{(2k)}}{(\phi'')^k}\,=\,1\cdot2\cdot\cdots (2k-1).$$
             In other words, the desired conclusion (\ref{normal}) follows by induction with respect to $k$. This concludes our proof of the asymptotic normality.

             \section{ Binary  Sequences not Equiprobable}\label{Binary}

             Consider again all $2^n$ sequences of length $n$ consisting merely of $0$'s and $1$'s. But assume now that the probability for a $0$ is $P(0) =p$ and the probability  for an $1$ is $P(1) =q.$ Here $p+q=1.$ For example, the sequence  $0\,1\,1\,0\,1\,1\,1$ has probability $p^2q^5$.
             The figurates in Table I and Table II (at the end) are constructed via (\ref{62}) below.
             The simple rule for the formation of these figurates is condensed in the formulae
             \begin{equation}\label{pqgen}
               f(x)\,=\,\begin{cases}
               (p+q)\underset {k=1}{\overset{\nu}{\prod}}\left(pqx^{-2k}+p^2+q^2+pqx^{2k}\right),\quad& n=2\nu+1\\
               \underset {k=1}{\overset{\nu}{\prod}}\left(pqx^{1-2k}+p^2+q^2+pqx^{2k-1}\right),    \quad& n=2\nu
             \end{cases}\end{equation}
             for the probability generating function. (Of course, the factor $p+q$ outside the product is $1$, but it is included to match Table II.)

             The characteristic function $\phi(\theta) = f(e^{i\theta})$ is
  \begin{equation}\label{pqchar}
               \phi(\theta)\,=\,\begin{cases}
               (p+q)\underset {k=1}{\overset{\nu}{\prod}}\left(p^2+q^2+2pq\cos(2k\theta)\right),\quad& n=2\nu+1\\
               \underset {k=1}{\overset{\nu}{\prod}}\left(p^2+q^2+2pq\cos((2k-1)\theta)\right),    \quad& n=2\nu
  \end{cases}\end{equation}
 For $p =q= 1/2$ we again obtain the expressions in Section \ref{Basic R}.

  Let us consider the case $n =2\nu+1$, the calculations for even $n$ being similar. Now
  \begin{equation*}  \phi'(\theta)\,=\,\phi(\theta)\,\sum_{k=1}^{\nu}\frac{-4pqk\sin(2k\theta)}{p^2+q^2+2pq\cos(2k\theta)}\,=\,\phi(\theta)A(\theta)
  \end{equation*}
  with an obvious abbreviation. The well-known expansion
  \begin{equation*}
    \frac{\varrho \sin(\psi)}{p^2+q^2-2\varrho\cos(\psi)}\,=\,\sum_{m=1}^{\infty}\varrho^m\sin(m\psi)\qquad (|\varrho| < 1)
  \end{equation*}
  converges for
  $$-\varrho\,=\,\min\{\frac{p}{q},\,\frac{q}{p}\}\quad\text{if}\quad p \neq q.$$
  Having treated the case $p=q=1/2$ in the previous sections, we assume that $p \neq q$ here. Then
   \begin{equation*}
    \frac{4k\varrho \sin(2k\theta)}{p^2+q^2-2\varrho\cos(2k\theta)}\,=\,4k\,\sum_{m=1}^{\infty}\varrho^m\sin(2km\theta)\qquad 
   \end{equation*}
   and so we obtain
   \begin{equation*} A(\theta)\,=\,  \sum_{m=1}^{\infty}\left(
     4\varrho^m  \sum_{k=1}^{\infty} k\sin(2km\theta)\right)\qquad (p\neq q).
   \end{equation*}
   Using the Maclaurin series for $\sin(2km\theta)$, we arrive at the formula
   \begin{equation}\label{sums} A(\theta)\,=\,\sum_{j=1}^{\infty}\left\{(-1)^{j+1}\frac{2^{2j+1}\bigl(1^{2j}+\dots+\nu^{2j}\bigr)}{(2j-1)!}\sum_{m=1}^{\infty}m^{2j-1}\varrho^m\right\}\theta^{2j-1}
   \end{equation}
   where some arrangements have been done. The corresponding convergence investigations are quite straightforward.

   By (\ref{sums}), $A(0)=0,\,A''(0) =0,\, A^{(4)}(0)=0,\dots,$ and
   \begin{equation}\label{touse}
     A^{(2j-1)}(0)\,=\,(-1)^{j+1}2^{2j+1}s_{\nu}(2j)\sum_{m=1}^{\infty}m^{2j-1}\varrho^m.
   \end{equation}
   (This expansion diverges for $\varrho = -1$, i.\,e. for $p=q$.)
   Here the infinite sum is easily calculated as the  differentiated geometric series  $$\sum_{m=1}^{\infty}m^{2j-1}\varrho^m\,=\,\Bigl(\varrho\frac{d}{\varrho}\Bigr)^{2j-1}\frac{1}{1-\varrho}\qquad (j=1,2,3,\dots).$$
   A calculation yields
   \begin{align*}
     A'(0)\,&=\,\frac{4\nu(\nu+1)(2\nu+1)}{3}\frac{\varrho}{(1-\varrho)^2},\\
     A'''(0)\,&=\, -\frac{16}{15} \nu(\nu+1)(2\nu+1)(3\nu^2+3\nu-1)\frac{1+4\varrho+\varrho^2}{(1-\varrho)^4}\varrho
   \end{align*}
   and using
   \begin{align*}
     \phi''(0)\,&=\,A'(0),\\
     \phi^{(4)}(0)\,&=\,3\phi''(0)A'(0)+A'''(0)\,=\,3A'(0)^2+A'''(0)
   \end{align*}
   we arrive at ({\ref{varpq}) and (\ref{fourpq}). ---The corresponding calculations for even $n$ yield the same final result.

   An analogous investigation as that in Section \ref{approach}, but now based on (\ref{touse}), shows the \emph{approach to normality} also for $p\neq q$. The difference is merely technical.

   \section{The Variance (with General $\ell$).}\label{The Vaiance}

   Consider again all sequences $x_1,x_2,\dots,x_n$ that can be formed of the digits $0,1,2,\dots,\ell$. Let
   $P_n(t;j_0,j_1,\dots,j_{\ell-1})$ count the number of those sequences  consisting of $j_0\,\,0$'s, $j_1\,\,1$'s,..., $j_{\ell-1}\,\,\ell$'s,\,\,$j_0+j_1+\dots+j_{\ell-1} =n$, for which $S=t$. For example
   $P_6(1;3,3) = 3,\,P_6(4;2,4)=2,\,P_9(27;3,3,3) =1$, and $P_9(t;3,3,3) = 0$ when $t \geq 28$.

   Constructing the sequence $x_1,x_2,\dots,x_n,x_{n+1}$ from $x_1,x_2,\dots,x_n$ we arrive at the fundamental recursive rule
   \begin{align}\label{recursion}
     P_{n+1}&(t;i_0,i_1,\dots,i_{\ell-1})\nonumber \\
     =&\,P_n(t-i_1-\dots-i_{\ell-1};i_0-1,i_1,\dots,i_{\ell-1})\nonumber\\
     &+\,P_n(t+i_0-i_2-\dots-i_{\ell-1};i_0,i_1-1,\dots,i_{\ell-1})+\cdots+\\
     &+\,P_n(t+i_0+\dots+i_{\ell-2};i_0,i_1,\dots,i_{\ell-2},i_{\ell-1}-1)\nonumber
   \end{align}
   where now $i_0+i_1+\dots+i_{\ell-1} = n+1.$ In passing, we notice that applying the recursive formula twice we obtain for $\ell=2$ that
   \begin{gather}\label{62}
     P_{n+1}(t;i_0,i_1)\,=\,P_{n-1}(t-2i_1;i_0-2,i_1)\nonumber\\
     +P_{n-1}(t+i_0-1;i_0-1,i_1-1) + P_{n-1}(t+i_0-i_1+1;i_0-1,i_1-1)\\
     +P_{n-1}(t+2i_0;i_0,i_1-2).\nonumber
   \end{gather}
   Repeated use of this identity describes exactly how the configurations in Tables I and II are built up.

   Let us return to the recursive rule (\ref{recursion}). For the calculation of the variance $\mu_2(n)$ we assume that all sequences  $x_1,x_2,\dots,x_n$ of length $n$ are equiprobable. Then
   $$\mu_2(n)\,=\,\sum_{x_{1}\cdots x_{n}}\frac{S^2(x_1,x_2,\dots ,x_n)}{\ell^n},$$
   where the sum is taken over all the $\ell^n$ possible sequences. The auxiliary quantities
   $$B_n(j_0,\dots,j_{\ell-1})\,=\,\sum_{t} t^2P_n(t;j_0,\dots,j_{\ell-1})$$
   satisfy according to the recursive rule (\ref{recursion}) the formula
   \begin{align*}
     B_{n+1}&(i_0,\dots,i_{\ell-1})\\
     =\,\sum_{k=0}^{\ell-1}&\sum_{t}\bigl(t+i_0+\dots+i_{k-1}-i_{k+1}-\dots-i_{\ell-1}\bigr)
     ^2\times\\ &P_n(t+i_0+\dots+i_{k-1}-i_{k+1}-\dots-i_{\ell-1};i_0,\dots,i_k-1,\dots,i_{\ell-1})\\
     -2\sum_{k=0}^{\ell-1}&\sum_{t}\bigl(i_0+\dots+i_{k-1}-i_{k+1}-\dots-i_{\ell-1}\bigr)\times\\
     &\bigl(t+i_0+\dots+i_{k-1}-i_{k+1}-\dots-i_{\ell-1}\bigr)\times\\
     &P_n(t+i_0+\dots+i_{k-1}-i_{k+1}-\dots-i_{\ell-1};i_0,\dots,i_k-1,\dots,i_{\ell-1})\\
     + \sum_{k=0}^{\ell-1}&\sum_{t}\bigl(i_0+i_1+\dots+i_{k-1}-i_{k+1}-\dots-i_{\ell-1}\bigr)
     ^2\times\\
  &P_n(t+i_0+\dots+i_{k-1}-i_{k+1}-\dots-i_{\ell-1};i_0,\dots,i_k-1,\dots,i_{\ell-1}),\\
   \end{align*}
   where $i_0+i_1+\dots+i_{\ell-1} = n+1.$ Here the first inner sum over $t$ is merely $B_n(i_0,i_1,\dots,i_{k-1},i_k-1,i_{k+1},\dots,i_{\ell-1})$ and the second inner sum is zero by symmetry. In the last inner sum $\sum_{t}P_n$ is a certain number of combinations.

   Therefore the above formula reduces to the simple expression
   \begin{align*}
     B_{n+1}&(i_0,\dots,i_{\ell-1})\,=\,\sum_{k=0}^{\ell-1}B_n(i_0,\dots,i_{k-1},i_k-1,i_{k+1},\dots,i_{\ell-1})\\
     &+\,\sum_{k=0}^{\ell-1}\bigl(i_0+\dots+i_{k-1}-i_{k+1}-\dots-i_{\ell-1}\bigr)^2\binom{n}{i_0\cdots i_{k-1}\,i_k-1\,i_{k+1}\cdots i_{\ell-1}}.
   \end{align*}

   Here
   $$\binom{n}{i_0\cdots i_{k}-1\cdots i_{\ell-1}}\,=\,\frac{n!}{i_0!\cdots (i_k-1)! \cdots i_{\ell-1}!}$$
   is the usual multinomial coefficient. Summing all equations with $i_0+i_1+\cdots+i_{\ell-1} =n+1$ and noting that
   $$\ell^n\mu_2(n)\,=\,\sum_{j_{0}+\cdots +j_{\ell-1}=n}B_n(j_0,\dots,j_{\ell-1}),$$
   we obtain
   \begin{align*}
     \mu_2(n+1)\,&=\,\mu_2(n)\\
    & +\,\ell^{-n-1}\sum_{k=0}^{\ell-1}\sum_{j_{0},\cdots ,j_{\ell-1}}\bigl(i_0+\cdots +i_{k-1}-i_{k+1}-\cdots-i_{l-1}\bigr)^2\times\\
     &\qquad\qquad\qquad\qquad \qquad\qquad\binom{n}{i_0\cdots i_{k}-1\cdots i_{\ell-1}}.
   \end{align*}
   The parenthesis in the sum can be written as the sum of products $\pm i_{\alpha}i_{\beta},$\,\,$\ell-1$ of which are of the form $i_{\alpha}^2,\,\,k^2+(\ell-k-1)^2$ of which are of the form $+i_{\alpha}i_{\beta}\,\,(\alpha\neq\beta)$, and $2k(\ell-k-1)$ of which are negative. Using well-known identities like
   \begin{align*}&\sum_{j_{0}+\cdots +j_{\ell-1}=n}j_0j_1\binom{n}{j_0\cdots j_{\ell-1}}\,=\,n(n-1)\ell^{n-2}\\
     &\sum_{j_{0}+\cdots +j_{\ell-1}=n}j_0^2\binom{n}{j_0\cdots j_{\ell-1}}\,=n\ell^{n-1}+n(n-1)\ell^{n-2}
   \end{align*}
   we finally obtain the equation
   \begin{equation}\label{addnow}
     \mu_2(n+1)\,=\,\mu_2(n) + n\,\frac{\ell-1}{\ell} +\frac{n(n-3)}{3}\,\frac{\ell^2-1}{\ell^2}.
   \end{equation}
   Adding the $n$ first equations (\ref{addnow}) and noting that $\mu_2(1) = 0$,
   we reach the final result (\ref{var}). This concludes our proof for the variance $\sigma^2(n)$.

   The fourth moment $\mu_4$ in (\ref{four}) is the result of a similar, although more tedious, calculation. However, a more effective method should be invented for higher moments. ---Corresponding formulae for $S^+$ are given in \cite{L1}.

   In passing we mention the formula
   $$S(i_0,\dots,i_{\ell-1})\,=\,2S^+(i_0,\dots,i_{\ell-1})-\frac{n^2-\bigl(i_0^2+\cdots+i_{\ell-1}^2\bigr)}{2},$$
   $i_0+\cdots +i_{\ell-1} =n$. Here the exceptional notation is understandable.

       \section{Edgeworth's Approximation}\label{Last}

       The closedness to normality of the distribution for $S$ is good, when $l^n$ is large. However, if $l^n$ is not large, especially the tails of the distribution behave obstinately, so that the assumption of normality is somewhat inadequate for precisely those values of $S$ whose significance may be in doubt. Fortunately, numerical calculations indicate that a correction based on Edgeworth's series gives an accurate approximation.

       Let
       $$\Phi(x)\,=\,\frac{1}{\sqrt{2\pi}}\int_{-\infty}^{x}e^{-\frac{t^2}{2}}\,dt$$
       denote the normal distribution function. The approximation
       \begin{equation}\label{edge}
         \mathrm{P}\{S\,\leq\,\sigma(n)x\}\,\approx\,\Phi(x) + \frac{1}{4!}\Bigl(\frac{\mu_4(n)}{\sigma(n)^4}-3\Bigr)\Phi^{(4)}(x)
       \end{equation}
       is obtained from Edgeworth's series [C, page 229], terms containing $\mu_6,\,\mu_7,\dots$ being neglected.

       It stands to reason that (\ref{edge}) is accurate, $\mu_4(n)$ and $\sigma(n)$ being calculated from ({\ref{four}) and (\ref{var}), provided that $\ell^n$ is large, say $\ell^n >10^6$. ---The dependence on $\ell$ is slightly puzzling. This is a point that requires further numerical investigation.
        \newpage
        {\centering
          \includegraphics[width=\textwidth]{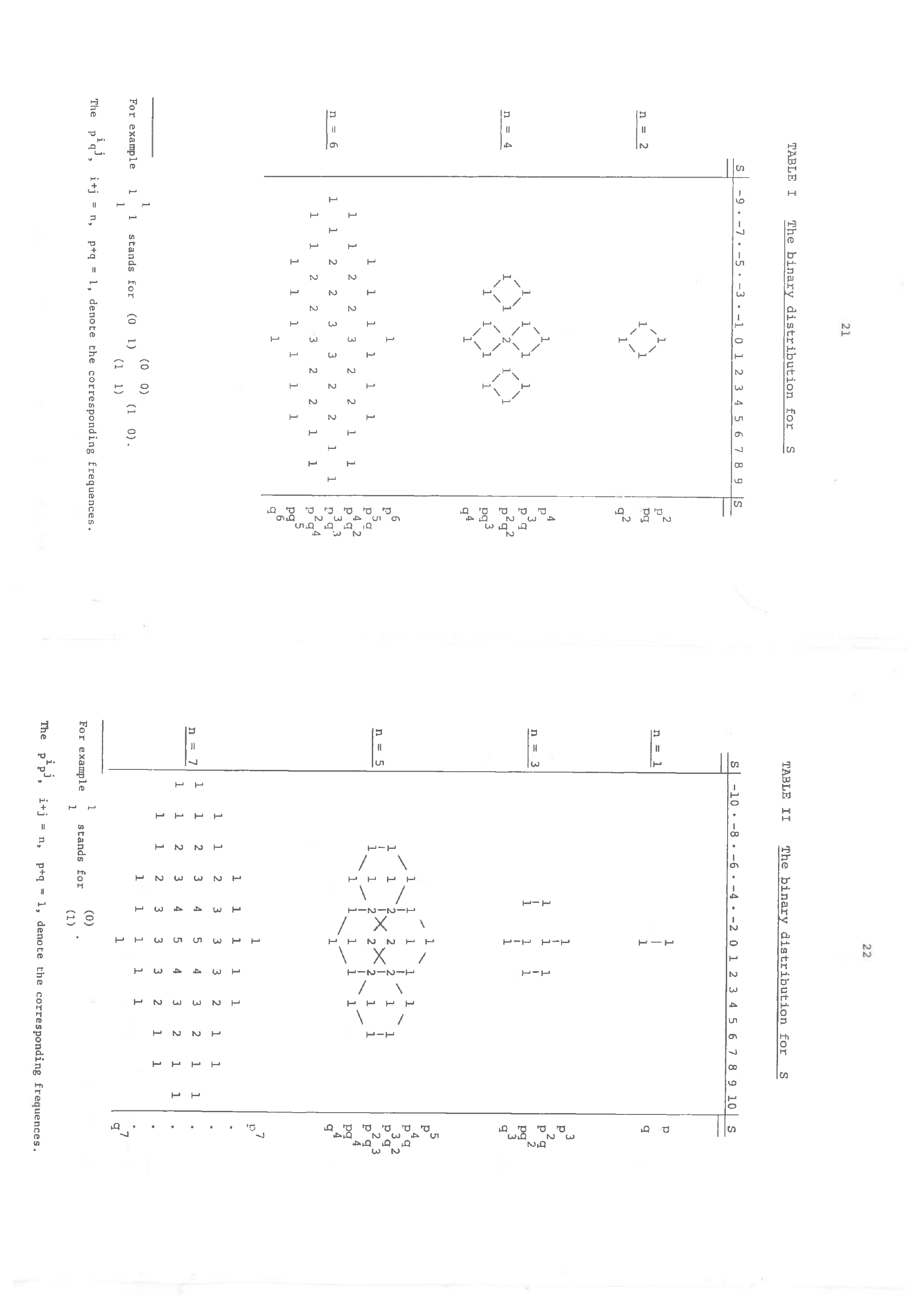}
           \par}
         
\end{document}